\documentclass[twosided,reqno]{amsart}

\usepackage{amsmath}
\usepackage{amssymb}
\usepackage{amsthm}

\theoremstyle{theorem}
\newtheorem{theorem}{\scshape Theorem }[section]

\theoremstyle{definition}

\numberwithin{equation}{section}

\begin{document}

\title{Identities involving Laguerre polynomials
derived from umbral calculus}

\author{Taekyun Kim}
\address{ Department of Mathematics, Kwangwoon University, Seoul 139-701, Republic of Korea.}
\email{tkkim@kw.ac.kr}


\maketitle

\begin{abstract}
In this paper, we investigate some identities of Laguerre polynomials involving Bernoulli and Euler polynomials which are derived from umbral calculus.
\end{abstract}

\section{Introduction}
For $\lambda(\neq 1) \in {\mathbb{C}}$, the {\it{Frobenius-Euler polynomials}} are given by the generating function to be
\begin{equation}\label{1}
\left(\frac{1-\lambda}{e^t-\lambda}\right)e^{xt}=e^{H(x|\lambda)t}=\sum_{n=0} ^{\infty} H_n (x|\lambda)\frac{t^n}{n!}{\text{ (see [6,7,13])}}
\end{equation}
with the usual convention about replacing $H^n (x|\lambda)$ by $H_n (x|\lambda)$.
In the special case, $x=0$, $H_n (0|\lambda)=H_n (\lambda)$ are called the {\it{$n$-th Frobenius-Euler numbers}}. By \eqref{1}, we easily see that
\begin{equation}\label{2}
H_n (x|\lambda)=\sum_{k=0} ^n \binom{n}{k}H_k(\lambda)x^{n-k},{\text{ (see [6])}}.
\end{equation}
As is well known, the {\it{Bernoulli polynomials}} are defined by the generating function to be
\begin{equation}\label{2}
\left(\frac{t}{e^t-1}\right)e^{xt}=e^{B(x)t}=\sum_{n=0} ^{\infty} B_n  (x) \frac{t^n}{n!},
\end{equation}
with the usual convention about replacing $B^n (x)$ by $B_n (x)$.
In the special case, $x=0$, $B_n (0)=B_n $ are called the {\it{$n$-th Bernoulli numbers}}. From \eqref{2}, we have
\begin{equation}\label{3}
B_n(x)=\sum_{k=0} ^n \binom{n}{k} B_k x^{n-k}.
\end{equation}
The {\it{Euler polynomials}} are also defined by the generating function to be
\begin{equation}\label{4}
\left(\frac{2}{e^t+1}\right)e^{xt}=\sum_{n=0} ^{\infty}E_n (x) \frac{t^n}{n!},{\text{ (see [2,3,8,10])}}.
\end{equation}
The Euler numbers, $E_n$, are defined by $E_n(0)=E_n$. By \eqref{4}, we get
\begin{equation}\label{5}
E_n(x)=\sum_{k=0} ^n \binom{n}{k} E_k x^{n-k},{\text{ (see [8,15,16,17])}}.
\end{equation}
As is well known, {\it{Laguerre polynomials}} are defined by the generating function to be
\begin{equation}\label{6}
\frac{\exp\left(-\frac{xt}{1-t}\right)}{1-t}=\sum_{n=0} ^{\infty}L_n(x)t^n,{\text{ (see[1,9])}}.
\end{equation}
From \eqref{6}, we note that
\begin{equation}\label{7}
\begin{split}
\sum_{n=0} ^{\infty} L_n(x)t^n &=\sum_{r=0} ^{\infty} \frac{(-1)^rx^r}{r!}(1-t)^{-r-1}t^r \\
&=\sum_{n=0} ^{\infty} \left(\sum_{r=0} ^n \frac{(-1)^r\binom{n}{r}}{r!}x^r\right)t^n.
\end{split}
\end{equation}
Thus, by \eqref{7}, we get
\begin{equation}\label{8}
L_n (x)=\sum_{r=0} ^n \frac{(-1)^r}{r!}\binom{n}{r}x^r.
\end{equation}
Thus, by \eqref{8}, we note that $L_n(x)$ is a polynomial of degree $n$ with rational coefficient and leading coefficient $\frac{(-1)^n}{n!}$. From \eqref{8}, we note that $u=L_n(x)$ is a solution of the following differential equation of order $2$:
\begin{equation}\label{9}
xu^{''}(x)+(1-x)u^{'}(x)+nu(x)=0.
\end{equation}
The {\it{Rodrigues' formula}} for $L_n(x)$ is given by
\begin{equation}\label{10}
L_n(x)=\frac{1}{n!}e^x \left(\frac{d^n}{dx^n}e^{-x}x^n\right),{\text{ (see [1,9])}}.
\end{equation}
By \eqref{10}, we easily see that
\begin{equation}\label{11}
\int_0 ^{\infty} e^{-x}L_m(x)L_n(x)dx=\delta_{m,n},~(m,n\in {\mathbb{N}}\cup\left\{0\right\}),
\end{equation}
where $\delta_{m,n}$ is the Kronecker symbol.

Let ${\mathbb{C}}$ be the complex number field and let ${\mathcal{F}}$ be the set of all formal power series in the variable $t$ over ${\mathbb{C}}$ with
\begin{equation*}
{\mathcal{F}}=\left\{ \left.f(t)=\sum_{k=0} ^{\infty} \frac{a_k}{k!} t^k~\right|~ a_k \in {\mathbb{C}} \right\}.
\end{equation*}
Let ${\mathbb{P}}={\mathbb{C}}[x]$. Then ${\mathbb{P}}^{*}$ will be defined by  the vector space of all linear functionals on ${\mathbb{P}}$. As is known, $\left< L~|~p(x)\right>$ denotes the action of a linear functional on a polynomial $p(x)$, and we remind that the vector space structure on ${\mathbb{P}}^{*}$ are defined by $\left<L+M|p(x) \right>=\left<L|p(x) \right>+\left<M|p(x) \right>$ and $\left<cL|p(x) \right>=c\left<L|p(x) \right>$, where $c$ is a complex constant. The formal power series
\begin{equation*}
f(t)=\sum_{k=0} ^{\infty} a_k \frac{t^k}{k!}\in {\mathcal{F}}
\end{equation*}
defines a linear functional on ${\mathbb{P}}$ by setting
\begin{equation}\label{12}
\left< f(t)|x^n\right>=a_n{\text{ for all }}n \geq 0{\text{ (see [4,5,7,13])}}.
\end{equation}
Thus, by \eqref{12}, we get
\begin{equation}\label{13}
\left<t^k | x^n \right>=n! \delta_{n,k},~(n,k \geq 0).
\end{equation}
Let $f_L(t)=\sum_{k=0} ^{\infty} \frac{\left<L|x^k\right>}{k!}t^k$. Then, by \eqref{13}, we see that $\left<f_L(t)|x^n\right>=\left<L|x^n\right>$ and so as a linear functional $L=f_L(t)$. Thus we note that the map $L \mapsto f_L(t)$ is a vector space isomorphism from ${\mathbb{P}}^{*}$ onto ${\mathcal{F}}$. Henceforth, ${\mathcal{F}}$ will denote both the algebra of formal power series in $t$ and the vector space of all linear functionals on ${\mathbb{P}}$, and so an element $f(t)$ of ${\mathcal{F}}$ will be thought of as both a formal power series and a linear functional. We shall call ${\mathcal{F}}$ the {\it{umbral algebra}}. The umbral calculus is the study of umbral algebra.

The order $o(f(t))$ of the non-zero power series $f(t)$ is the smallest integer $k$ for which the coefficient of $t^k$ does not vanish (see \cite{13,14}). If $o(f(t))=1$, then $f(t)$ is called a {\it{delta series}} and if $o(f(t))=0$, then $f(t)$ is called an {\it{invertible series}}. From \eqref{12} and \eqref{13}, we note that $\left<e^{yt}|x^n\right>=y^n$ and so $\left<e^{yt}|p(x)\right>=p(y)$. Let $f(t)\in{\mathcal{F}}$ and $p(x)\in{\mathbb{P}}$. Then we have the following equations (see \cite{11,12,13}):
\begin{equation}\label{14}
f(t)=\sum_{k=0} ^{\infty} \frac{\left<f(t)|x^k\right>}{k!} t^k,~p(x)=\sum_{k=0} ^{\infty} \frac{\left<t^k|p(x)\right>}{k!}x^k.
\end{equation}
It is easy to show that
\begin{equation*}
\left<f(t)g(t)|p(x)\right>=\left<f(t)|g(t)p(x)\right>=\left<g(t)|f(t)p(x)\right>,
\end{equation*}
where $f(t),g(t)\in {\mathcal{F}}$ and $p(x)\in{\mathbb{P}}$. From \eqref{14}, we can derive the following equations:
\begin{equation}\label{15}
p^{(k)}(0)=\left<t^k|p(x)\right>{\text{ and }}\left<t^0|p^{(k)}(x)\right>=p^{(k)}(0).
\end{equation}
Thus, by \eqref{15}, we get
\begin{equation}\label{16}
t^kp(x)=p^{(k)}(x)=\frac{d^kp(x)}{dx^k},~e^{yt}p(x)=p(x+y).
\end{equation}
Let $o(f(t))=1$ and $o(g(t))=0$. Then there exists a unique sequence $S_n(x)$
such that
\begin{equation}\label{17}
\left<g(t)f(t)^k|S_n(x)\right>=n!\delta_{n,k}, ~(n,k \geq 0).
\end{equation}
This sequence $S_n(x)$ is called the {\it{Sheffer sequence}} for $(g(t),f(t))$, which is denoted by $S_n(x)\sim (g(t),f(t))$. Let $S_n(x)\sim(g(t),f(t))$. Then we have the following equations (see \cite{06,07,13}) :
\begin{equation}\label{18}
h(t)=\sum_{k=0} ^{\infty} \frac{\left<h(t)|S_k(x)\right>}{k!}g(t)f(t)^k,~p(x)=\sum_{k=0} ^{\infty} \frac{\left<g(t)f(t)^k|p(x)\right>}{k!}S_k(x),
\end{equation}
\begin{equation}\label{19}
f(t)S_n(x)=nS_{n-1}(x),~\left<f(t)|p(\alpha x)\right>=\left<f(\alpha t)|p(x)\right>,
\end{equation}
and
\begin{equation}\label{20}
\frac{1}{g({\bar{f}}(t))}e^{y{\bar{f}}(t)}=\sum_{k=0} ^{\infty} \frac{S_k(y)}{k!}t^k,{\text{ for all }}y \in {\mathbb{C}},
\end{equation}
where ${\bar{f}}(t)$ is the compositional inverse of $f(t)$.

In this paper, we investigate some identities of Laguerre polynomials involving Bernoulli and Euler polynomials which are derived from umbral calculus.

\section{Some identities of Laguerre polynomials}

Let
\begin{equation}\label{21}
{\mathbb{P}}_n=\left\{p(x)\in {\mathbb{C}}[x]|\deg p(x) \leq n \right\}.
\end{equation}
Then we note that ${\mathbb{P}}_n$ is an inner product space with weighted inner product
\begin{equation}\label{22}
\left<p(x),q(x)\right>=\int_0 ^{\infty} e^{-x}p(x)q(x)dx,
\end{equation}
where $p(x),q(x)\in{\mathbb{P}}_n$. By \eqref{11}, \eqref{21} and \eqref{22}, we easily wee that $L_0(x), L_1(x),\ldots,L_n(x)$ are orthogonal basis for ${\mathbb{P}}_n$.

For $p(x)\in{\mathbb{P}}_n$, let us assume that polynomial $p(x)$ defined on $[0,\infty)$ is given by
\begin{equation}\label{23}
p(x)=\sum_{k=0} ^{\infty} C_kL_k(x).
\end{equation}
Then, by \eqref{11} and \eqref{23}, we get
\begin{equation}\label{24}
\begin{split}
C_k&=\left<p(x),L_k(x)\right>=\int_0 ^{\infty} e^{-x}p(x)L_k(x)dx \\
&=\frac{1}{k!}\int_0 ^{\infty} \left(\frac{d^k}{dx^k}e^{-x}x^k\right)p(x)dx.
\end{split}
\end{equation}
Let us take $p(x)=x^n$ $(n\geq 0)$. Then, by \eqref{24}, we get
\begin{equation}\label{25}
\begin{split}
C_k&=\frac{1}{k!}\int_0 ^{\infty}\left(\frac{d^k}{dx^k}e^{-x}x^k\right)x^ndx \\
&=(-1)^k\frac{(n)_k}{k!}\int_0 ^{\infty} e^{-x}x^ndx=(-1)^k\binom{n}{k}n!.
\end{split}
\end{equation}
From \eqref{23} and \eqref{25}, we have
\begin{equation}\label{26}
x^n=n!\sum_{k=0} ^n (-1)^k\binom{n}{k}L_k(x).
\end{equation}
Note that
\begin{equation}\label{27}
(x)_k=x(x-1)\cdots(n-k+1)\sim\left(1,e^t-1\right).
\end{equation}
By \eqref{18} and \eqref{26}, we easily get
\begin{equation}\label{28} x^n=\sum_{k=0} ^n \frac{\left.\left<(e^t-1)^k\right|x^n\right>}{k!}(x)_k=\sum_{k=0} ^n S_2(n,k)(x)_k,
\end{equation}
where $S_2(n,k)$ is the Stirling number of the second kind.

From \eqref{26} and \eqref{28}, we have
\begin{equation*}
L_k(x)=\frac{(-1)^k}{n!\binom{n}{k}}S_2(n,k)(x)_k,
\end{equation*}
where $n,k \in {\mathbb{Z}}_+$ with $n \geq k$.
\begin{theorem}\label{thm1}
For $n,k \in {\mathbb{Z}}_+$ with $n \geq k$, we have
\begin{equation*}
L_k(x)=\frac{(-1)^k}{n!\binom{n}{k}}S_2(n,k)(x)_k.
\end{equation*}
\end{theorem}
For $p_n(x)$ and $q_n(x)=\sum_{k=0} ^n q_{n,k}x^k$ sequences of polynomials, we define the umbral composition of $q_n(x)$ with $p_n(x)$ to be the sequences
\begin{equation}\label{29}
\left(q_n \circ p\right)(x)=\sum_{k=0} ^n q_{n,k} p_k(x),{\text{ (see [3,6,7,13])}}.
\end{equation}
From \eqref{6} and \eqref{20}, we note that
\begin{equation}\label{30}
L_n(x)\sim\left(1-t,\frac{t}{1-t}\right).
\end{equation}
Let $s_n(x)\sim(g(t),f(t))$, $r_n(x)\sim(h(t),l(t))$. Then we have
\begin{equation}\label{31}
\left(r_n \circ s \right)(x)=(g(t)h(f(t)),l(f(t))),{\text{ (see [13])}}.
\end{equation}
From \eqref{30} and \eqref{31}, we can derive
\begin{equation}\label{32}
\left(L_n \circ L\right)(x)\sim(1,t).
\end{equation}
By \eqref{20} and \eqref{32}, we get
\begin{equation}\label{33}
\left(L_n \circ L\right)(x)=x^n.
\end{equation}
Therefore, by \eqref{26} and \eqref{33}, we obtain the following theorem.
\begin{theorem}\label{thm2}
For $n \geq 0$, we have
\begin{equation*}
\frac{1}{n!}\left(L_n\circ L\right)(x)=\sum_{k=0} ^n (-1)^k\binom{n}{k}L_k(x).
\end{equation*}
\end{theorem}
For $p(x)=E_n(x)\in {\mathbb{P}}_n$, let us assume that
\begin{equation}\label{34}
E_n(x)=\sum_{k=0} ^n C_k L_k (x).
\end{equation}
By \eqref{24} and \eqref{34}, we get
\begin{equation}\label{35}
\begin{split}
C_k&=\frac{1}{k!}\int_0 ^{\infty} \left(\frac{d^k}{dx^k}e^{-x}x^k \right)E_n(x)dx \\
&=(-1)^k\binom{n}{k} \int_0 ^{\infty}e^{-x}x^k E_{n-k} (x) dx \\
&=(-1)^k \binom{n}{k} \sum_{l=0} ^{n-k} \binom{n-k}{l}\int_0 ^{\infty}e^{-x}x^{k+l}dx \\
&=n!\sum_{l=0} ^{n-k} (-1)^k \binom{k+l}{l}\frac{E_{n-k-l}}{(n-k-l)!}.
\end{split}
\end{equation}
From \eqref{4}, we have
\begin{equation}\label{36}
\left(\frac{2}{e^t+1}\right)=\left(\frac{e^t+1}{2}\right)^{-1}=\left(1+\frac{e^t-1}{2}\right)^{-1}=\sum_{j=0} ^{\infty}(-1)^j\left(\frac{e^t-1}{2}\right)^j.
\end{equation}
Note that $(e^t-1)^j$ is a delta series and $E_n(x)\sim \left(\frac{e^t+1}{2},t\right)$. By \eqref{4}, \eqref{16} and \eqref{20}, we get
\begin{equation}\label{37}
\begin{split}
\sum_{n=0} ^{\infty}E_n(x)\frac{t^n}{n!}&=\left(\frac{2}{e^t+1}\right)e^{xt}=\sum_{j=0} ^{\infty} (-1)^j\left(\frac{e^t-1}{2}\right)^je^{xt} \\
&=\sum_{n=0} ^{\infty}\left(\sum_{j=0} ^n (-1)^j \left(\frac{e^t-1}{2}\right)^jx^n\right)\frac{t^n}{n!} \\
&=\sum_{n=0} ^{\infty}\left(\sum_{j=0} ^n \left(-\frac{1}{2}\right)^j(e^t-1)^jx^n \right)\frac{t^n}{n!}.
\end{split}
\end{equation}
By comparing the coefficients on the both sides of \eqref{37}, we get
\begin{equation}\label{38}
\begin{split}
E_n(x)&=\sum_{j=0} ^n \left(-\frac{1}{2}\right)^j(e^t-1)^jx^n \\
&=\sum_{j=0} ^n \left(-\frac{1}{2}\right)^jj!\sum_{k=j} ^{\infty} S_2(k,j)\frac{t^k}{k!}x^n \\
&=\sum_{j=0} ^n \left(-\frac{1}{2}\right)^jj!\sum_{k=j} ^{\infty} S_2(k,j)\binom{n}{k}x^{n-k} \\
&=\sum_{k=0} ^n \binom{n}{k}\left\{\sum_{j=0} ^k \left(-\frac{1}{2}\right)^jj!S_2(k,j)\right\}x^{n-k}.
\end{split}
\end{equation}
From \eqref{5} and \eqref{38}, we have
\begin{equation}\label{39}
E_k=\sum_{j=0} ^k\left(-\frac{1}{2}\right)^jj!S_2(k,j).
\end{equation}
By \eqref{35} and \eqref{39}, we get
\begin{equation}\label{40}
\begin{split}
C_k&=n!\sum_{l=0} ^{n-k} (-1)^k\binom{k+l}{l}\frac{1}{(n-k-l)!}\sum_{j=0} ^{n-k-l} \left(-\frac{1}{2}\right)^jj!S_2(n-k-l,j) \\
&=n!\sum_{l=0} ^{n-k}\sum_{j=0} ^{n-k-l}(-1)^{k+l}\binom{k+l}{l}\frac{j!S_2(n-k-l,j)}{2^j(n-k-l)!}.
\end{split}
\end{equation}
Therefore, by \eqref{34} and \eqref{40}, we obtain the following theorem.
\begin{theorem}\label{thm3}
For $n \geq 0$, we have
\begin{equation*}
E_n(x)=n!\sum_{k=0} ^n\sum_{l=0} ^{n-k}\sum_{j=0} ^{n-k-l}(-1)^{k+l}\binom{k+l}{l}\frac{j!S_2(n-k-l,j)}{2^j(n-k-l)!}L_k(x).
\end{equation*}
\end{theorem}
By \eqref{2}, we see that
\begin{equation}\label{36}
\left(\frac{t}{e^t-1}\right)=\left(\frac{e^t-1}{t}\right)^{-1}=\left(1+\frac{e^t-t-1}{t}\right)^{-1}=\sum_{j=0} ^{\infty}(-1)^j\left(\frac{e^t-t-1}{t}\right)^j.
\end{equation}
Note that $\left(\frac{e^t-t-1}{t}\right)^j=\left(\frac{e^t-1}{t}-1\right)^j$ is a delta series in ${\mathcal{F}}$ and $B_n(x)\sim\left(\frac{e^t-1}{t},t\right)$. By \eqref{2}, \eqref{20} and \eqref{36}, we get
\begin{equation}\label{42}
\begin{split}
\sum_{n=0} ^{\infty}B_n(x)\frac{t^n}{n!}&=\left(1+\frac{e^t-t-1}{t}\right)^{-1}e^{xt}=\sum_{j=0} ^{\infty}(-1)^j\left(\frac{e^t-t-1}{t}\right)^je^{xt} \\
&=\sum_{n=0} ^{\infty}\left(\sum_{j=0} ^n (-1)^j\left(\frac{e^t-t-1}{t}\right)^j x^n \right)\frac{t^n}{n!}.
\end{split}
\end{equation}
By comparing the coefficients on the both sides of \eqref{42}, we get
\begin{equation}\label{43}
B_n(x)=\sum_{j=0} ^n (-1)^j \left(\frac{e^t-t-1}{t}\right)^j x^n.
\end{equation}
From \eqref{14}, we have
\begin{equation}\label{44}
\begin{split}
\left(\frac{e^t-t-1}{t}\right)^j x^n&=\sum_{k=0} ^{n-j} \frac{\left<t^k \left| \left(\frac{e^t-t-1}{t}\right)^j x^n\right.\right>}{k!}x^k \\
&=\sum_{k=0} ^{n-j}\frac{\left.\left< \left(\frac{e^t-t-1}{t}\right)^j \right| t^kx^n\right>}{k!}x^k \\
&=\sum_{k=0} ^{n-j} \binom{n}{k}\left.\left< \left(\frac{e^t-t-1}{t}\right)^j \right| x^{n-k}\right>x^k \\
&=\sum_{k=0} ^{n-j} \binom{n}{k}\sum_{l=0} ^j \binom{j}{l}(-1)^{j-l} \left<t^0\left| \left(\frac{e^t-1}{t}\right)^lx^{n-k} \right.\right>x^k.
\end{split}
\end{equation}
It is easy to show that
\begin{equation}\label{45}
\begin{split}
\left(\frac{e^t-1}{t}\right)^lx^{n-k}&=\frac{1}{t^l}l!\sum_{m=l} ^{\infty}S_2(m,l)\frac{t^m}{m!}x^{n-k} \\
&=\sum_{m=0} ^{n-k}S_2(m+l,l)\frac{l!}{(m+l)!}t^mx^{n-k} \\
&=\sum_{m-0} ^{n-k} S_2(m+l,l)\frac{l!}{(m+l)!}(n-k)_mx^{n-k-m}.
\end{split}
\end{equation}
Thus, by \eqref{44} and \eqref{45}, we get
\begin{equation}\label{46}
\left(\frac{e^t-t-1}{t}\right)^jx^n=\sum_{k=0} ^{n-j} \sum_{l=0} ^j \binom{n}{k}\binom{j}{l}(-1)^{j-l}\frac{S_2(n-k+l,l)}{\binom{n-k+l}{l}}x^k.
\end{equation}
Thus, by \eqref{43} and \eqref{46}, we get
\begin{equation}\label{47}
\begin{split}
B_n(x)&=\sum_{j=0} ^n \sum_{k=0} ^{n-j} \sum_{l=0} ^j \binom{n}{k}\binom{j}{l}(-1)^l\frac{S_2(n-k+l,l)}{\binom{n-k+l}{l}}x^k \\
&=\sum_{k=0} ^n \binom{n}{k} \left\{\sum_{j=0} ^{n-k} \sum_{l=0} ^j\binom{j}{l}(-1)^l\frac{S_2(n-k+l,l)}{\binom{n-k+l}{l}}\right\}x^k\\
&=\sum_{k=0} ^n \binom{n}{k} \left\{\sum_{j=0} ^k \sum_{l=0} ^j (-1)^l\binom{j}{l}\frac{S_2(k+l,l)}{\binom{k+l}{l}}(-1)^l\right\}x^{n-k}.
\end{split}
\end{equation}
By \eqref{3} and \eqref{47}, we get
\begin{equation}\label{48}
B_k=\sum_{j=0} ^k \sum_{l=0} ^j (-1)^l\binom{j}{l}\frac{S_2(k+l,l)}{\binom{k+l}{l}}.
\end{equation}
From \eqref{3}, we note that $B_n(x)\in{\mathbb{P}}_n$. So, let us assume that
\begin{equation}\label{49}
B_n(x)=\sum_{k=0} ^n C_kL_k(x).
\end{equation}
By \eqref{24}, we get
\begin{equation}\label{50}
\begin{split}
C_k&=\frac{1}{k!}\int_0 ^{\infty} \left(\frac{d^k}{dx^k}e^{-x}x^k \right)B_n(x)dx \\
&=(-1)^k\binom{n}{k} \int_0 ^{\infty}e^{-x}x^k B_{n-k} (x) dx \\
&=(-1)^k \binom{n}{k} \sum_{l=0} ^{n-k} \binom{n-k}{l}B_{n-k-l}\int_0 ^{\infty}e^{-x}x^{k+l}dx \\
&=(-1)^k \binom{n}{k}\sum_{l=0} ^{n-k} \binom{n-k}{l}B_{n-k-l} (k+l)! \\
&=n!\sum_{l=0} ^{n-k} \frac{(-1)^kB_{n-k-l}}{(n-k-l)!} \binom{k+l}{l}.
\end{split}
\end{equation}
From \eqref{48} and \eqref{50}, we have
\begin{equation}\label{51}
C_k=n!\sum_{l=0} ^{n-k}\sum_{j=0} ^{n-k-l} \sum_{m=0} ^j \frac{(-1)^{m+k}\binom{j}{m}S_2(n-k-l+m,m)}{(n-k-l)!\binom{n-k-l+m}{m}}.
\end{equation}
Therefore, by \eqref{49} and \eqref{51}, we obtain the following theorem.
\begin{theorem}\label{thm4}
For $n \geq 0$, we have
\begin{equation*}
B_n(x)=n!\sum_{k=0} ^n \left\{\sum_{l=0} ^{n-k}\sum_{j=0} ^{n-k-l} \sum_{m=0} ^j \frac{(-1)^{m+k}\binom{j}{m}S_2(n-k-l+m,m)}{(n-k-l)!\binom{n-k-l+m}{m}}\right\}L_k(x).
\end{equation*}
\end{theorem}
From \eqref{1}, we note that
\begin{equation}\label{52}
\frac{1-\lambda}{e^t-\lambda}=\left(1+\frac{e^t-1}{1-\lambda}\right)^{-1}=\sum_{j=0} ^{\infty}(-1)^j\left(\frac{1}{1-\lambda}\right)^j(e^t-1)^j.
\end{equation}
Note that $\left(e^t-1\right)^j$ is a delta series in ${\mathcal{F}}$ and $H_n(x|\lambda)\sim\left(\frac{e^t-\lambda}{1-\lambda},t\right)$.
By \eqref{1}, \eqref{20} and \eqref{52}, we get
\begin{equation}\label{53}
H_n(x|\lambda)=\sum_{j=0} ^n (-1)^j\left(\frac{1}{1-\lambda}\right)^j(e^t-1)^jx^n,
\end{equation}
and
\begin{equation}\label{54}
\begin{split}
(e^t-1)^jx^n&=j!\sum_{k=j} ^n S_2(k,j)\frac{t^k}{k!}x^n \\
&=j!\sum_{k=j} ^n S_2(k,j)\binom{n}{k}x^{n-k}.
\end{split}
\end{equation}
From \eqref{53} and \eqref{54}, we can derive the following equation \eqref{55}:
\begin{equation}\label{55}
\begin{split}
H_n(x|\lambda)&=\sum_{j=0} ^n \sum_{k=j} ^n \binom{n}{k}(-1)^j\frac{j!}{(1-\lambda)^j}S_2(k,j)x^{n-k} \\
&=\sum_{k=0} ^n \binom{n}{k}\left\{\sum_{j=0} ^k\frac{j!}{(\lambda-1)^j}S_2(k,j)\right\}x^{n-k},
\end{split}
\end{equation}
and
\begin{equation}\label{56}
H_n(x|\lambda)=\sum_{k=0} ^n \binom{n}{k} H_k(\lambda)x^{n-k}.
\end{equation}
thus, by \eqref{55} and \eqref{56}, we get
\begin{equation}\label{57}
H_k(\lambda)=\sum_{j=0} ^k \frac{j!}{(\lambda-1)^j}S_2(k,j).
\end{equation}
By \eqref{56}, we see that $H_n(x|\lambda)\in {\mathbb{P}}_n$. Let us assume that
\begin{equation}\label{58}
H_n(x|\lambda)=\sum_{k=0} ^n C_k L_k(x).
\end{equation}
From \eqref{24}, we have
\begin{equation}\label{59}
\begin{split}
C_k&=\frac{1}{k!}\int_0 ^{\infty} \left(\frac{d^k}{dx^k}e^{-x}x^k \right)H_n(x|\lambda)dx \\
&=(-1)^k\binom{n}{k} \int_0 ^{\infty}e^{-x}x^k H_{n-k} (x|\lambda) dx \\
&=(-1)^k \binom{n}{k} \sum_{l=0} ^{n-k} \binom{n-k}{l}H_{n-k-l}(\lambda)\int_0 ^{\infty}e^{-x}x^{k+l}dx \\
&=(-1)^k \binom{n}{k}\sum_{l=0} ^{n-k} \binom{n-k}{l}H_{n-k-l}(\lambda) (k+l)! \\
&=n!\sum_{l=0} ^{n-k} \binom{k+l}{l}(-1)^k\frac{H_{n-k-l}(\lambda)}{(n-k-l)!}.
\end{split}
\end{equation}
By \eqref{57} and \eqref{59}, we get
\begin{equation}\label{60}
C_k=n!\sum_{l=0} ^{n-k}\sum_{j=0} ^{n-k-l} \frac{(-1)^k\binom{k+l}{l}j!}{(n-k-l)!(\lambda-1)^j}S_2(n-k-l,j).
\end{equation}
Therefore, by \eqref{58} and \eqref{60}, we obtain the following theorem.
\begin{theorem}\label{thm5}
For $n \geq 0$, we have
\begin{equation*}
H_n(x|\lambda)=n!\sum_{k=0} ^n \left\{\sum_{l=0} ^{n-k}\sum_{j=0} ^{n-k-l}\binom{k+l}{l}(-1)^k\frac{j!S_2(n-k-l,j)}{(n-k-l)!(\lambda-1)^j}\right\}L_k(x).
\end{equation*}
\end{theorem}
For $s_n(x)\sim(g(t),f(t))$ and $r_n(x)\sim(h(t),l(t))$, let us assume that
\begin{equation}\label{61}
s_n(x)=\sum_{k=0} ^n C_{n,k} r_k(x).
\end{equation}
Then, we note that
\begin{equation}\label{62}
C_{n,k}=\frac{1}{k!}\left. \left<\frac{h({\bar{f}}(t))}{g({\bar{f}}(t))}l({\bar{f}}(t))^k \right|x^n \right>
\end{equation}
where $n,k \geq 0$ (see \cite{06,07,13}). From \eqref{4}, \eqref{6} and \eqref{20}, we not that
\begin{equation}\label{63}
L_n(x)\sim\left(1-t,\frac{t}{t-1}\right),~E_n(x)\sim\left(\frac{1+e^t}{2},t\right).
\end{equation}
Now, we assume that
\begin{equation}\label{64}
E_n(x)=\sum_{k=0} ^n C_{n,k} L_k(x).
\end{equation}
By \eqref{62} and \eqref{63}, we get
\begin{equation}\label{65}
\begin{split}
C_{n,k}&=-\frac{1}{k!}\left.\left<\left(\frac{2}{e^t+1}\right)\frac{t^k}{(t-1)^{k-1}}\right|x^n \right> \\
&=-\binom{n}{k}\left.\left<\frac{2}{e^t+1}\right|\left(\frac{1}{t-1}\right)^{k-1}x^{n-k}\right>.
\end{split}
\end{equation}
It is easy to show that
\begin{equation}\label{66}
\begin{split}
\left(\frac{1}{t-1}\right)^{k-1}x^{n-k}&=\sum_{l=0} ^{\infty}\binom{k-l-2}{l}(-1)^{k-1}t^lx^{n-k} \\
&=\sum_{l=0} ^{n-k} \binom{k-l-2}{l}(-1)^{k-1}(n-k)_lx^{n-k-l}.
\end{split}
\end{equation}
By \eqref{65} and \eqref{66}, we get
\begin{equation}\label{67}
\begin{split}
C_{n,k}&=\binom{n}{k}(-1)^k\sum_{l=0} ^{n-k}\binom{k-l-2}{l}\frac{(n-k)!}{(n-k-l)!}\left.\left<\frac{2}{e^t+1}\right|x^{n-k-l}\right> \\
&=n!\sum_{l=0} ^{n-k} \binom{k-l-2}{l}(-1)^k\frac{1}{k!(n-k-l)!}E_{n-k-l}.
\end{split}
\end{equation}
Therefore, by \eqref{64} and \eqref{67}, we obtain the following theorem.
\begin{theorem}\label{thm6}
For $n \geq 0$, we have
\begin{equation*}
E_n(x)=n!\sum_{k=0} ^n \left\{\sum_{l=0} ^{n-k} \binom{k-l-2}{l}(-1)^k\frac{E_{n-k-l}}{k!(n-k-l)!}\right\}L_k(x).
\end{equation*}
\end{theorem}
By Theorem \ref{thm3} and Theorem \ref{thm6}, we get
\begin{equation*}
\sum_{l=0} ^{n-k} \binom{k-l-2}{l}\frac{E_{n-k-l}}{k!(n-k-l)!}=\sum_{l=0} ^{n-k}\sum_{j=0} ^{n-k-l} (-1)^l\binom{k+l}{l}\frac{j!S_2(n-k-l,j)}{2^j(n-k-l)!},
\end{equation*}
where $n,k \geq 0$ with $n \geq k$. Let us assume that
\begin{equation}\label{68}
B_n(x)=\sum_{k=0} ^n C_{n,k} L_k(x).
\end{equation}
By \eqref{62}, we get
\begin{equation}\label{69}
\begin{split}
C_{n,k}&=-\frac{1}{k!}\left.\left<\left(\frac{e^t-1}{t}\right)^{-1}\frac{t^k}{(t-1)^{k-1}}\right|x^n \right> \\
&=-\frac{1}{k!}\left.\left<\left(\frac{t}{e^t-1}\right)\frac{1}{(t-1)^{k-1}}\right|t^kx^n \right> \\
&=-\binom{n}{k}\sum_{l=0} ^{n-k} \binom{k-l-2}{l}(-1)^{k-1}\frac{(n-k)!}{(n-k-l)!}\left.\left<\frac{t}{e^t-1}\right|x^{n-k-l}\right>
\\
&=n!\sum_{l=0} ^{n-k} \binom{k-l-2}{l}(-1)^k\frac{1}{k!(n-k-l)!}B_{n-k-l}.
\end{split}
\end{equation}
Therefore, by \eqref{68} and \eqref{69}, we obtain the following theorem.
\begin{theorem}\label{thm7}
For $n \geq 0$, we have
\begin{equation*}
B_n(x)=n!\sum_{k=0} ^n \left\{\sum_{l=0} ^{n-k} \binom{k-l-2}{l}(-1)^k\frac{B_{n-k-l}}{k!(n-k-l)!}\right\}L_k(x).
\end{equation*}
\end{theorem}
Now, we set
\begin{equation}\label{70}
H_n(x|\lambda)=\sum_{k=0} ^n C_{n,k} L_k(x).
\end{equation}
By \eqref{62}, we get
\begin{equation}\label{71}
\begin{split}
C_{n,k}&=-\frac{1}{k!}\left.\left<\left(\frac{1-\lambda}{e^t-\lambda}\right)\frac{t^k}{(t-1)^{k-1}}\right|x^n\right> \\
&=n!\sum_{l=0} ^{n-k} \binom{k-l-2}{l}(-1)^k\frac{H_{n-k-l} (\lambda)}{k!(n-k-l)!}.
\end{split}
\end{equation}
From \eqref{70} and \eqref{71}, we have
\begin{equation*}
H_n(x|\lambda)=n!\sum_{k=0} ^n \left\{\sum_{l=0} ^{n-k} \binom{k-l-2}{l}(-1)^k\frac{H_{n-k-l} (\lambda)}{k!(n-k-l)!}\right\}L_k(x).
\end{equation*}

\end{document}